\documentclass[10pt]{amsart}
\usepackage{amsmath}
\usepackage{amsxtra}
\usepackage{amscd}
\usepackage{amsthm}
\usepackage{amsfonts}
\usepackage{amssymb}
\usepackage{eucal}
\usepackage[matrix,arrow,curve]{xy}

\newtheorem{thm}{Theorem}[section]
\newtheorem{cor}[thm]{Corollary}
\newtheorem{lem}[thm]{Lemma}
\newtheorem{prop}[thm]{Proposition}
\newtheorem{conj}[thm]{Conjecture}
\theoremstyle{remark}
\newtheorem{rem}[thm]{Remark}

\theoremstyle{definition}

\theoremstyle{remark}

\newcommand{\nc}{\newcommand}
\nc{\renc}{\renewcommand} \nc{\ssec}{\subsection}
\nc{\sssec}{\subsubsection} \nc{\on}{\operatorname}

\nc\ol{\overline} \nc\ul{\underline} \nc\wt{\widetilde}
\nc\tboxtimes{\wt{\boxtimes}} \nc{\alp}{\alpha}

\nc{\ZZ}{{\mathbb Z}} \nc{\NN}{{\mathbb N}} \nc{\CC}{{\mathbb C}}
\nc{\OO}{{\mathbb O}} \renc{\SS}{{\mathbb S}} \nc{\DD}{{\mathbb
D}}
\nc{\Fq}{{\mathbb F}_q} \nc{\Fqb}{\ol{{\mathbb F}_q}}
\nc{\Ql}{\ol{{\mathbb Q}_\ell}} \nc{\id}{\text{id}} \nc\X{\mathcal
X}

\nc{\Hom}{\on{Hom}} \nc{\Lie}{\on{Lie}} \nc{\Loc}{\on{Loc}}
\nc{\Pic}{\on{Pic}} \nc{\Bun}{\on{Bun}} \nc{\IC}{\on{IC}}
\nc{\Aut}{\on{Aut}} \nc{\rk}{\on{rk}} \nc{\Sh}{\on{Sh}}
\nc{\Perv}{\on{Perv}} \nc{\pos}{{\on{pos}}} \nc{\Conv}{\on{Conv}}
\nc{\Sph}{\on{Sph}} \nc{\Sym}{\on{Sym}}
\nc{\BunBb}{\overline{\Bun}_B} \nc{\Buno}{\overset{o}{\Bun}}
\nc{\BunPb}{{\overline{\Bun}_P}}
\nc{\BunBM}{\overline{\Bun}_{B(M)}}
\nc{\BunPbw}{{\widetilde{\Bun}_P}}
\nc{\BunBP}{\widetilde{\Bun}_{B,P}} \nc{\GUb}{\overline{G/U}}
\nc{\GUPb}{\overline{G/U(P)}}
\nc{\iso}{{\stackrel{\sim}{\longrightarrow}}}

\nc{\Hhom}{\underline{\on{Hom}}} \nc\syminfty{\on{Sym}^{\infty}}
\nc\lal{\ol{\lambda}} \nc\xl{\ol{x}} \nc\thl{\ol{\theta}}
\nc\nul{\ol{\nu}} \nc\mul{\ol{\mu}} \nc\Sum\Sigma
\nc{\oX}{\overset{o}{X}{}}

\nc{\M}{{\mathcal M}} \nc{\N}{{\mathcal N}} \nc{\F}{{\mathcal F}}
\nc{\D}{{\mathcal D}} \nc{\Q}{{\mathcal Q}} \nc{\Y}{{\mathcal Y}}
\nc{\G}{{\mathcal G}} \nc{\E}{{\mathcal E}} \nc{\CalC}{{\mathcal
C}}
\nc\Dh{\widehat{\D}}

\nc{\C}{{\mathcal C}} \nc{\K}{{\mathcal K}}
\renewcommand{\H}{{\mathcal H}}

\nc{\T}{{\mathcal T}} \nc{\V}{{\mathcal V}} \renc{\P}{{\mathcal
P}} \nc{\A}{{\mathcal A}} \nc{\B}{{\mathcal B}} \nc{\U}{{\mathcal
U}}

\nc{\Gr}{\on{Gr}}

\nc{\frn}{{\check{\mathfrak u}(P)}}
\nc\f{{\mathfrak f}}

\nc{\q}{{\mathfrak q}} \nc{\p}{{\mathfrak p}} \nc{\s}{{\mathfrak
s}} \nc\w{\text{w}}

\nc\mathi\iota \nc\Spec{\on{Spec}} \nc\Mod{\on{Mod}}
\nc{\tw}{\widetilde{\mathfrak t}} \nc{\pw}{\widetilde{\mathfrak
p}} \nc{\qw}{\widetilde{\mathfrak q}} \nc{\jw}{\widetilde j}

\nc{\grb}{\overline{\Gr}} \nc{\I}{\mathcal I}

\nc{\lambdach}{{\check\lambda}} \nc{\Lambdach}{{\check\Lambda}{}}
\nc{\much}{{\check\mu}} \nc{\omegach}{{\check\omega}}
\nc{\nuch}{{\check\nu}} \nc{\etach}{{\check\eta}}
\nc{\alphach}{{\check\alpha}} \nc{\betach}{{\check\beta}}
\nc{\rhoch}{{\check\rho}} \nc{\ch}{{\check h}}

\nc{\Hb}{\overline{\H}}

\nc{\BA}{{\mathbb{A}}} \nc{\BC}{{\mathbb{C}}}
\nc{\BG}{{\mathbb{G}}} \nc{\BQ}{{\mathbb{Q}}}
\nc{\BM}{{\mathbb{M}}} \nc{\BN}{{\mathbb{N}}} \nc{\BO}{{\mathbb{O}}}
\nc{\BP}{{\mathbb{P}}} \nc{\BR}{{\mathbb{R}}}
\nc{\BZ}{{\mathbb{Z}}} \nc{\BS}{{\mathbb{S}}}

\nc{\CA}{{\mathcal{A}}} \nc{\CB}{{\mathcal{B}}}
\nc{\CE}{{\mathcal{E}}} \nc{\CF}{{\mathcal{F}}}
\nc{\CG}{{\mathcal{G}}} \nc{\CH}{{\mathcal{H}}}
\nc{\CI}{{\mathcal{I}}} \nc{\CL}{{\mathcal{L}}}
\nc{\CM}{{\mathcal{M}}} \nc{\CN}{{\mathcal{N}}}
\nc{\CO}{{\mathcal{O}}} \nc{\CP}{{\mathcal{P}}}
\nc{\CQ}{{\mathcal{Q}}} \nc{\CR}{{\mathcal{R}}}
\nc{\CS}{{\mathcal{S}}} \nc{\CT}{{\mathcal{T}}}
\nc{\CU}{{\mathcal{U}}} \nc{\CV}{{\mathcal{V}}}  \nc{\CY}{{\mathcal Y}}
\nc{\CW}{{\mathcal{W}}} \nc{\CZ}{{\mathcal{Z}}}

\nc{\cM}{{\check{\mathcal M}}{}} \nc{\csM}{{\check{\mathcal A}}{}}
\nc{\oM}{{\overset{\circ}{\mathcal M}}{}}
\nc{\obM}{{\overset{\circ}{\mathbf M}}{}}
\nc{\oCA}{{\overset{\circ}{\mathcal A}}{}}
\nc{\obA}{{\overset{\circ}{\mathbf A}}{}}
\nc{\ooM}{{\overset{\circ}{M}}{}}
\nc{\osM}{{\overset{\circ}{\mathsf M}}{}}
\nc{\vM}{{\overset{\bullet}{\mathcal M}}{}}
\nc{\nM}{{\underset{\bullet}{\mathcal M}}{}}
\nc{\oD}{{\overset{\circ}{\mathcal D}}{}}
\nc{\obD}{{\overset{\circ}{\mathbf D}}{}}
\nc{\oA}{{\overset{\circ}{\mathbb A}}{}}
\nc{\op}{{\overset{\bullet}{\mathbf p}}{}}
\nc{\cp}{{\overset{\circ}{\mathbf p}}{}}
\nc{\oU}{{\overset{\bullet}{\mathcal U}}{}}
\nc{\oZ}{{\overset{\circ}{\mathcal Z}}{}}
\nc{\ofZ}{{\overset{\circ}{\mathfrak Z}}{}}

\nc{\ff}{{\mathfrak{f}}} \nc{\fv}{{\mathfrak{v}}}
\nc{\fa}{{\mathfrak{a}}} \nc{\fb}{{\mathfrak{b}}}
\nc{\fd}{{\mathfrak{d}}} \nc{\fe}{{\mathfrak{e}}}
\nc{\fg}{{\mathfrak{g}}} \nc{\fgl}{{\mathfrak{gl}}}
\nc{\fh}{{\mathfrak{h}}} \nc{\fri}{{\mathfrak{i}}}
\nc{\fj}{{\mathfrak{j}}} \nc{\fk}{{\mathfrak{k}}}
\nc{\fm}{{\mathfrak{m}}} \nc{\fn}{{\mathfrak{n}}}
\nc{\ft}{{\mathfrak{t}}} \nc{\fu}{{\mathfrak{u}}}
\nc{\fw}{{\mathfrak{w}}} \nc{\fz}{{\mathfrak{z}}}
\nc{\fp}{{\mathfrak{p}}} \nc{\fq}{{\mathfrak{q}}}
\nc{\frr}{{\mathfrak{r}}}
\nc{\fs}{{\mathfrak{s}}} \nc{\fo}{{\mathfrak{o}}}
\nc{\fsl}{{\mathfrak{sl}}} \nc{\fsp}{{\mathfrak{sp}}}
\nc{\hsl}{{\widehat{\mathfrak{sl}}}}
\nc{\hgl}{{\widehat{\mathfrak{gl}}}}
\nc{\hg}{{\widehat{\mathfrak{g}}}}
\nc{\chg}{{\widehat{\mathfrak{g}}}{}^\vee}
\nc{\hn}{{\widehat{\mathfrak{n}}}}
\nc{\chn}{{\widehat{\mathfrak{n}}}{}^\vee}

\nc{\fA}{{\mathfrak{A}}} \nc{\fB}{{\mathfrak{B}}} \nc{\fC}{{\mathfrak{C}}}
\nc{\fD}{{\mathfrak{D}}} \nc{\fE}{{\mathfrak{E}}}
\nc{\fF}{{\mathfrak{F}}} \nc{\fG}{{\mathfrak{G}}} \nc{\fH}{{\mathfrak{H}}}
\nc{\fI}{{\mathfrak{I}}} \nc{\fJ}{{\mathfrak{J}}}
\nc{\fK}{{\mathfrak{K}}} \nc{\fL}{{\mathfrak{L}}}
\nc{\fM}{{\mathfrak{M}}} \nc{\fN}{{\mathfrak{N}}}
\nc{\fP}{{\mathfrak{P}}} \nc{\fQ}{{\mathfrak{Q}}}
\nc{\fT}{{\mathfrak{T}}} \nc{\fU}{{\mathfrak{U}}}
\nc{\fV}{{\mathfrak{V}}} \nc{\fW}{{\mathfrak{W}}}
\nc{\fX}{{\mathfrak{X}}} \nc{\fY}{{\mathfrak{Y}}}
\nc{\fZ}{{\mathfrak{Z}}}

\nc{\ba}{{\mathbf{a}}}
\nc{\bb}{{\mathbf{b}}} \nc{\bc}{{\mathbf{c}}}
\nc{\be}{{\mathbf{e}}} \nc{\bj}{{\mathbf{j}}}
\nc{\bn}{{\mathbf{n}}} \nc{\bp}{{\mathbf{p}}}
\nc{\bq}{{\mathbf{q}}} \nc{\br}{{\mathbf{r}}} \nc{\bt}{{\mathbf{t}}}
\nc{\bfu}{{\mathbf{u}}} \nc{\bv}{{\mathbf{v}}}
\nc{\bx}{{\mathbf{x}}} \nc{\by}{{\mathbf{y}}}
\nc{\bw}{{\mathbf{w}}} \nc{\bA}{{\mathbf{A}}}
\nc{\bB}{{\mathbf{B}}} \nc{\bC}{{\mathbf{C}}}
\nc{\bD}{{\mathbf{D}}} \nc{\bF}{{\mathbf{F}}}
\nc{\bH}{{\mathbf{H}}} \nc{\bK}{{\mathbf{K}}}
\nc{\bM}{{\mathbf{M}}} \nc{\bN}{{\mathbf{N}}}
\nc{\bO}{{\mathbf{O}}} \nc{\bS}{{\mathbf{S}}} \nc{\bT}{{\mathbf{T}}}
\nc{\bV}{{\mathbf{V}}} \nc{\bW}{{\mathbf{W}}}
\nc{\bX}{{\mathbf{X}}}
\nc{\bY}{{\mathbf{Y}}} \nc{\bP}{{\mathbf{P}}}
\nc{\bZ}{{\mathbf{Z}}} \nc{\bh}{{\mathbf{h}}}
\nc{\bmu}{{\boldsymbol{\mu}}}

\nc{\sA}{{\mathsf{A}}} \nc{\sB}{{\mathsf{B}}}
\nc{\sC}{{\mathsf{C}}} \nc{\sD}{{\mathsf{D}}}
\nc{\sE}{{\mathsf{E}}} \nc{\sF}{{\mathsf{F}}} \nc{\sG}{{\mathsf{G}}}
\nc{\sK}{{\mathsf{K}}} \nc{\sL}{{\mathsf{L}}}
\nc{\sM}{{\mathsf{M}}} \nc{\sO}{{\mathsf{O}}}
\nc{\sQ}{{\mathsf{Q}}} \nc{\sP}{{\mathsf{P}}}
\nc{\sT}{{\mathsf{T}}} \nc{\sZ}{{\mathsf{Z}}}
\nc{\sV}{{\mathsf{V}}}
\nc{\sfp}{{\mathsf{p}}} \nc{\sr}{{\mathsf{r}}}
\nc{\st}{{\mathsf{t}}} \nc{\sfb}{{\mathsf{b}}}
\nc{\sfc}{{\mathsf{c}}} \nc{\sd}{{\mathsf{d}}}
\nc{\sz}{{\mathsf{z}}}

\nc{\BK}{{\bar{K}}}

\nc{\tA}{{\widetilde{\mathbf{A}}}}
\nc{\tB}{{\widetilde{\mathcal{B}}}}
\nc{\tg}{{\widetilde{\mathfrak{g}}}} \nc{\tG}{{\widetilde{G}}}
\nc{\TM}{{\widetilde{\mathbb{M}}}{}}
\nc{\tO}{{\widetilde{\mathsf{O}}}{}}
\nc{\tU}{{\widetilde{\mathfrak{U}}}{}} \nc{\TZ}{{\tilde{Z}}}
\nc{\tx}{{\tilde{x}}} \nc{\tbv}{{\tilde{\bv}}}
\nc{\tfP}{{\widetilde{\mathfrak{P}}}{}} \nc{\tz}{{\tilde{\zeta}}}
\nc{\tmu}{{\tilde{\mu}}}

\nc{\urho}{\underline{\rho}} \nc{\uB}{\underline{B}}
\nc{\uC}{{\underline{\mathbb{C}}}} \nc{\ui}{\underline{i}}
\nc{\uj}{\underline{j}} \nc{\ofP}{{\overline{\mathfrak{P}}}}
\nc{\oB}{{\overline{\mathcal{B}}}}
\nc{\og}{{\overline{\mathfrak{g}}}} \nc{\oI}{{\overline{I}}}

\nc{\eps}{\varepsilon} \nc{\hrho}{{\hat{\rho}}}
\nc{\blambda}{{\boldsymbol{\lambda}}}

\nc{\one}{{\mathbf{1}}} \nc{\two}{{\mathbf{t}}}

\nc{\Rep}{{\mathop{\operatorname{\rm Rep}}}}
\nc{\Tot}{{\mathop{\operatorname{\rm Tot}}}}
\nc{\Ker}{{\mathop{\operatorname{\rm Ker}}}}
\nc{\Hilb}{{\mathop{\operatorname{\rm Hilb}}}}
\nc{\End}{{\mathop{\operatorname{\rm End}}}}
\nc{\Ext}{{\mathop{\operatorname{\rm Ext}}}}
\nc{\CHom}{{\mathop{\operatorname{{\mathcal{H}}\it om}}}}
\nc{\GL}{{\mathop{\operatorname{\rm GL}}}}
\nc{\gr}{{\mathop{\operatorname{\rm gr}}}}
\nc{\Id}{{\mathop{\operatorname{\rm Id}}}}
\nc{\defi}{{\mathop{\operatorname{\rm def}}}}
\nc{\length}{{\mathop{\operatorname{\rm length}}}}
\nc{\supp}{{\mathop{\operatorname{\rm supp}}}}

\nc{\Cliff}{{\mathsf{Cliff}}}
\nc{\Fl}{{\mathsf{Fl}}} \nc{\Fib}{{\mathsf{Fib}}}
\nc{\Coh}{{\mathsf{Coh}}} \nc{\FCoh}{{\mathsf{FCoh}}}

\nc{\reg}{{\text{\rm reg}}}

\nc{\cplus}{{\mathbf{C}_+}} \nc{\cminus}{{\mathbf{C}_-}}
\nc{\cthree}{{\mathbf{C}_*}} \nc{\Qbar}{{\bar{Q}}}

\nc{\bOmega}{{\overline{\Omega}}}

\nc{\seq}[1]{\stackrel{#1}{\sim}}

\nc{\aff}{\operatorname{aff}}

\begin{document}

\author{Michael Finkelberg and Dmitry Kubrak}
\title
%[Dva Idiota]
{Vanishing cycles on Poisson varieties}

\dedicatory{To Vadim Schechtman on his 60th birthday}

%\thanks{{\bf Mathematics Subject Classification (2000).}
%19E08, (22E65, 37K10).}

%\thanks{{\bf Key words.} $q$-difference Toda lattice, Equivariant
%$K$-theory, Laumon compactification.}

%\thanks{The work of L.R. was partially supported
%by  RFBR grants 07-01-92214-CNRSL-a and 05-01-02805-CNRSL-a.
%L.R. gratefully acknowledges the support from Deligne 2004 Balzan
%prize in mathematics.}

\address{{\it Address}:\newline
National Research University Higher School of Economics, \newline
Department of Mathematics; \newline
20 Myasnitskaya st,
Moscow 101000, Russia; \newline
Institute for Information Transmission Problems of RAS}

\email{\newline fnklberg@gmail.com
\newline dmkubrak@gmail.com}

\begin{abstract}
We extend slightly the results of Evens-Mirkovi\'c, and ``compute''
the characteristic cycles of Intersection Cohomology sheaves on the
transversal slices in the double affine Grassmannian.
We propose a conjecture relating the hyperbolic stalks
and the microlocalization at a torus-fixed point in a Poisson variety.
\end{abstract}
\maketitle

%\tableofcontents

\section{Introduction}

\subsection{}
For a perverse sheaf $\CF$ on a smooth complex variety $Y$, a basic problem
in microlocal geometry is to compute the characteristic cycle
$\on{CC}(\CF)=\sum_{S\in\CS}c_S\ol{T^*_SY}$: a positive integral combination
of the closures of the conormal bundles to locally closed strata $S\subset Y$.
This problem was solved by S.~Evens and I.~Mirkovi\'c in~\cite{em} for
Goresky-MacPherson sheaves of $G[[z]]$-orbits in the affine Grassmannian
$\on{Gr}_G$ of a semisimple complex group $G$. Namely, they were able to
prove the vanishing of all the Euler obstructions by computing the
$\sA$-equivariant Euler characteristic of the complex links ($\sA\subset G$
is a Cartan torus). In their approach, it was essential that the strata
$S$ were the orbits of some algebraic group action, and that the conormal
fibers could be identified with some algebraic Lie algebras.

The results of~\cite{em} suggest the existence of a microlocal fiber functor
from the geometric Satake category $\on{Perv}_{G[[z]]}(\on{Gr}_G)$ to the
category of representations of the Langlands dual group $\check G$.
A conjectural affine analogue $\on{Gr}^{\on{aff}}_G$ of $\on{Gr}_G$
(the double affine Grassmannian) was introduced in~\cite{bf}. More precisely,
the whole $\on{Gr}^{\on{aff}}_G$ is too infinite-dimensional to handle by the
machinery available at the moment, and we are bound to think of it as of a
collection of (finite-dimensional) transversal slices. For one thing, the
price we have to pay for this lame approach is that the strata are no longer
orbits of any group action.

Fortunately, they are still the symplectic leaves of a Poisson structure.
This equips the conormal fibers with a Lie algebra structure, and the
argument of Evens-Mirkovi\'c goes through. This is proved
in~Section~\ref{a} after some general preparation in~Section~\ref{Euler},
answering a question of S.~Evens.

An isomorphism of the (would be) microlocal fiber functor with the
standard Mirkovi\'c-Vilonen fiber functor is a relation between the
hyperbolic stalks and the microlocalization at the torus-fixed points.
It seems to go through in our extended setup, see~Conjecture~\ref{oko}.
If true, it would express the classical $R$-matrices of~\cite[Section~4.8]{mo}
as the action of the microlocal fundamental group.

There is one more rich class of algebraic Poisson varieties with finitely
many symplectic leaves: Nakajima quiver varieties. They satisfy all the
technical conditions of~\ref{setup} except one (apart from type $A$ case,
finite or affine): the torus fixed point set is not discrete. So already
in the simplest examples of type $D,E$ Kleinian singularities (finite type
quiver varieties, corresponding to the zero weight of the adjoint fundamental
representation) the Euler obstruction does not vanish. More precisely, it
is equal to 1 for all type $D,E$ Kleinian singularities (there are but 2
strata), so that the multiplicity of the cotangent fiber at the singular
point in the characteristic cycle of IC (=constant) sheaf is 2. On the positive
side, for affine hypertoric varieties (another class of algebraic Poisson varieties with finitely many symplectic leaves) all our conditions are
satisfied, so that the multiplicity of the conormal bundle to a stratum in
the characteristic cycle of the IC sheaf coincides with the total dimension of
the IC stalk at this stratum, computed in~\cite{pw}.

\subsection{Acknowledgments}
This note is a result of patient explanations by V.~Ginzburg, A.~Kuznetsov,
L.~Rybnikov, P.~Schapira, G.~Williamson.
The research of M.F. was carried out at the IITP RAS at the expense of the 
Russian Foundation for Sciences (project no. 14-50-00150).
The article was prepared by D.K. within the framework of a subsidy 
granted to the HSE by the Government of the Russian Federation for the 
implementation of the Global Competitiveness Program.

%\documentclass[12 pt, engish, a4paper]{article}
%\usepackage{amsfonts, amsmath, latexsym, amssymb,euscript}
%\newcounter{vse}[section]
%\newcommand*{\theo}{\par\noindent\addtocounter{vse}{1}%
%\textbf {Theorem \arabic{section}.\arabic{vse}.} }

%\newcommand*{\prop}{\par\noindent\addtocounter{vse}{1}%
%\textbf {Proposition \arabic{section}.\arabic{vse}.} }

%\newcommand*{\lem}{\par\noindent\addtocounter{vse}{1}%
%\textbf {Lemma \arabic{section}.\arabic{vse}.} }

%\newcommand*{\cor}{\par\noindent\addtocounter{vse}{1}%
%\textbf {Corollary \arabic{section}.\arabic{vse}.} }

\newcommand*{\mf}{\mathfrak}
\newcommand*{\mb}[1]{\mathbb{#1}}
\newcommand*{\mr}{\mathrm}
\newcommand*{\mc}{\mathcal}
\newcommand*{\g}[1]{\mathfrak{#1}}
\newcommand*{\Oo}{\mathcal O}
\newcommand*{\eu}{\EuScript}

%Let $M$ be an affine Poisson variety and $\mc S= \{S_i\}$ ---
%stratification of $M$ by symplectic leaves. Let also one of the leaves be
%a point $p$ and and assume that there is a $\mb C^\times$-action contracting
% $M$ to $p$. We consider a situation when there is an additional action of
%reductive group G, preserving the Poisson bracket and such that there is a
%subtori $T\hookrightarrow G$ such that $p$ is isolated fixed point of
%$T$-compact form $T^c$. We prove then that the Euler characteristic of
%complex link of $p$ and any other leaf $S_i$ is equal to zero. Applying
%this to quiver varieties we prove vanishing of all Euler obstructions.

%Let $\mc F$ be a perverse sheaf on a complex manifold $X$, constructible
%with respect to a Whitney stratification $\mc L$. A basic problem in
%microlocal geometry is to compute the characteristic cycle of $\mc F$:
%$$
%CC(\mc F)=\sum_{\alpha\in \mc L} c_{\alpha}(\mc F)\cdot
%\overline{T_{\alpha}^*(X)}
%$$

%The main reason for studying complex links is Dubson's and Kashiwara's index
%formula

\section{Vanishing of the Euler obstruction}
\label{Euler}

\subsection{Review of~\cite{em}}
\label{index}
We recall the argument of S.~Evens and I.~Mirkovi\'c.

Let $\EuScript L$ be a Whitney stratification of a complex manifold $Y$ and let $D_{\EuScript L}(Y)$ denote the derived category of $\EuScript L$-constructible sheaves. Let $\EuScript F$ be an object of $D_\CL(Y)$. Then $\EuScript F$ defines two important constructible (with respect to $\EuScript L$) functions on $Y$. The first one, $\chi$ is just the Euler characteristic $\chi(\eu F_y)$ of the stalk of $\eu F$ at a point $y$. As $\eu F$ lies in $D_{\EuScript L}(Y)$, $\chi$ is $\eu L$-constructible and we denote by $\chi_{\alpha}(\eu F)$ its value on a stratum $\alpha$. The second function, $c(\eu F)$ comes from the characteristic cycle
$$
CC(\eu F)=\sum_{\alpha\in \eu L} c_{\alpha}(\eu F)\cdot \overline{T_{\alpha}^*(Y)},
$$
its value on stratum $\alpha$ is equal to $c_{\alpha}(\eu F)$, which is called
the microlocal multiplicity of $\eu F$ along $\alpha$. Here $T_{\alpha}^* (Y)\subset T^*(Y)$ denotes the conormal bundle to the stratum $\alpha$.

  For any pair of strata $\alpha$ and $\beta$ the Euler obstruction $e_{\alpha,\beta}$ is defined as $c_{\alpha}(\mb C_{\beta})$, where $\mb C_\beta$ is the constant sheaf on $\beta$, extended by zero on $Y$.
%As we will see $e_{\alpha,\beta}$ are rather clear topological
%invariants of $Y$.

  A covector $\xi \in T_{\alpha}^*(Y)$ is called $\eu L$-generic if it lies in

  $$
  T_{\alpha}^*(Y)^r:=T_{\alpha}^*(Y) - \cup_{\alpha\neq\beta}\overline{T_{\beta}^*(Y)}.
  $$
The set of generic elements is open and dense in $T_{\alpha}^*(Y)$ and its fundamental group is called microlocal fundamental group of the stratum $\alpha$.

For the following theorem, see e.g.~\cite[Chapter~9]{ks}:
\begin{thm}
\label{dubson}
a) On the Grothendieck group of $\eu L$-constructible complexes,
$$
c_\alpha = \sum_{\beta\in\eu L} e_{\alpha,\beta}\cdot \chi_{\beta}
$$

b) One has $e_{\alpha,\alpha}=(-1)^{\mr{dim}(\alpha)}$ and $e_{\alpha,\beta}=0$ if $\alpha\nsubseteq\overline{\beta}$. If $\alpha\subset\partial\overline{\beta}$, we choose
\begin{itemize}
\item a normal slice $N$ in $(Y,\eu L)$ to $\alpha$ at a point $y\in\alpha$,
\item a holomorphic function $\phi$ on $N$ vanishing at $y$ and such that
$d_y\phi\in T_y^*(N)$ is $\eu L$-generic,
\item a small ball $B$ around $y$ in $Y$,
\item a small $t\in \mb C$.
\end{itemize}
Then
$$
e_{\alpha,\beta}=(-1)^{\mr{dim}(\alpha)+1} \chi_c (\beta \cap \phi^{-1} (t)\cap B)
$$
is, up to sign, the Euler characteristic of a compactly supported cohomology of the intersection of $\beta$ with a nearby hyperplane $\phi^{-1}(t)$, near
$\alpha$ and normal to $\alpha$.
\end{thm}

The intersection $\beta \cap \phi^{-1} (t)\cap B$ is called the complex link of strata $\alpha$ and $\beta$.

Now Mirkovi\'c and Evens argue as follows: having a compact torus ${\sA^c}$ acting on $Y$ stabilizing $y$ and $\beta$,  choose $B$ in~Theorem~\ref{dubson}
${\sA^c}$-invariant and try to find a ${\sA^c}$-invariant function $\phi$ on $N$ such that its differential is generic. Then the complex link will be stable under ${\sA^c}$-action too and the classical result of Borel can be applied:
Let $Z$ be a paracompact space with finite cohomological dimension with the action of a compact torus ${\sA^c}$. Then
$\chi_c(Z)=\chi_c(Z^{\sA^c})$, where $Z^{\sA^c}$ is the fixed-point set.
In particular if $Z^{\sA^c}$ is empty, we get $\chi_c(Z)=0$.
In other words, $e_{\alpha,\beta}=0$.

%This is actually the reason why  we have asked for $p$ being the
%\textit{isolated} fixed point of some torus-action. Then the link
%between $p$ and any other symplectic leaf will be endowed with the torus
%action without fixed points and so its euler characteristic will be zero.

So the main question left is the existence of an ${\sA^c}$-invariant function
with generic differential. In the next subsection we formulate some
technical conditions which guarantee the existence. In the next section we
check the technical conditions in our main example.

\subsection{Setup}
\label{setup}
$X$ is an affine complex algebraic Poisson variety with finitely many
symplectic leaves $X_i,\ i\in L$.
To avoid a misunderstanding, let us comment
on the definition of symplectic leaves in this situation. Let
$X_{\on{sing}}\subset X$ be the singular locus of $X$. Then $X-X_{\on{sing}}$
is a smooth Poisson variety, and we require it to have finitely many
symplectic leaves. Furthermore, $X_{\on{sing}}$ inherits a Poisson structure
from $X$, but $\dim X_{\on{sing}}<\dim X$, so the symplectic leaves are defined
on $X_{\on{sing}}$ (by induction in the dimension). Note that the symplectic
leaves are the locally closed algebraic subvarieties of $X$.
The closure of $X_i$ is denoted by $\ol{X}_i$.
The Poisson bracket on $R:=\BC[X]$ is denoted by $\{,\}$.
A reductive group $G$ with Lie algebra $\fg$ acts on $X$
preserving the bracket $\{,\}$.
We assume the existence of a moment map $\mu:\ X\to\fg^*$, hence
$\mu^*:\ \Sym(\fg)=\BC[\fg^*]\to R$.
The fixed point subset $X^\sA$ for a maximal torus $\sA\subset G$ consists of
a single point $x\in X$, a 0-dimensional symplectic leaf.
Additionally, a one-parametric group $\BC^\times$ acts on $X$ commuting with
$G$. We assume that the corresponding action on $R$ has only positive weights:
it induces the grading $R=\BC\oplus\bigoplus_{l>0}R_l$,
and $\fm_x=\bigoplus_{l>0}R_l$. We assume the Poisson bracket has weight
$w>0$ with respect to the one-parametric group $\BC^\times$, i.e.
$\{R_l,R_m\}\subset R_{l+m-w}$. Finally, we assume that $R_l=0$ for $0<l<w$,
and $R_w$ coincides with the image of $\mu^*$ on linear functions on $\fg^*:\
\fg\stackrel{\mu^*}{\to}R$.

\subsection{Cotangent Lie algebra}
\label{Cot}
%Let $M$ be an affine Poisson variety over $\mb C$ and let $A=\mb C[M]$ be
%algebra of functions on $M$. Let $\{f,g\}$ denote the Poisson bracket of
%$f,g\in A$ and let $\mc S=\{S_i\} $ be the stratification of $M=\sqcup S_i$
%by the symplectic leaves. We consider the situation, when one of the leaves
%is just a point $p$.
%It is equivalent to the fact that Poisson bivector defining the bracket is
%equal to zero in $p$. We also need an additional $\mb C^{\times}$-action
%$\varphi$ on $M$, that rescales the Poisson bracket and contracts $M$ to $p$.
%Then $\varphi$ induces a natural grading on $A$ by the means of
%$\mb C^\times$-action:
% $$
% A=\bigoplus_{n\in \mb Z} A_n, \mbox{ and }\varphi(t)|A_n = t^n A_n
% $$
%This grading is obviously multiplicative and if Poisson bracket is rescaled
%by $t^k$ then $\{A_i,A_j\}\subset A_{i+j-k}$. In our case, as $ \varphi $ is
%contracting whole $M$ to a point, we have $A_i=0$ for all $i<0$ and
%$A_0=\mb C$.
The Poisson bracket on $R$ induces a Lie algebra structure on the Zariski
cotangent space $T^*_xX:=\fm_x/\fm_x^2$.
%We claim that in this situation the Zariski cotangent space $T^*_{Zar,p}$ to
%$M$ in point $p$ has a natural graded Lie algebra structure, coming from
%Poisson bracket on $A$. The Zariski cotangent space is by definition equal
%to $\mf m_p/\mf m_p^2$, where $\mf m_p$ is maximal ideal of $A$, corresponding
%to $p$. In our case $\mf m_p$ is naturally identified with
%$A_{>0} =\oplus_{i>0} A_i$.
In effect, since $x$ is a symplectic leaf, the Poisson bivector vanishes at
$x$, and $\{\fm_x,\fm_x\}\subset\fm_x$. By Leibniz rule
$\{fg,h\} =f\{g,h\}+g\{f,h\}$, we also have
$\{\fm_x^2,\fm_x\}\subset \fm_x^2$. In other words, the Poisson bracket
on $\fm_x$ descends to $\fm_x/\fm_x^2$ and provides it with a structure of
Lie algebra, to be denoted $\fq$. The grading on $\fm_x$ gives rise to a
grading on $\fq=\bigoplus_{l\geq w}\fq_w$. We define $\fq^m:=\fq_{w+m}$.
With this new shifted grading, the Lie bracket is homogeneous, so $\fq$
becomes a nonnegatively graded Lie algebra.

\begin{lem}
\label{alg}
$\fq$ is an algebraic Lie algebra, i.e. there exists a complex algebraic
Lie group $\sQ$ with Lie algebra $\fq$.
\end{lem}

\proof
By our assumptions in~\ref{setup}, $\fq^0$ is a quotient of $\fg$,
moreover, it is an algebraic quotient, i.e. $\fq^0=\on{Lie}\sQ^0$ for a
Lie group quotient $G\twoheadrightarrow\sQ^0$. In effect,
let $\fk:=\on{Ker}\mu^*\subset\fg$. Then $\xi\in\fg$ lies in $\fk$ iff
it acts by the zero vector field on the Zariski tangent cone
$C_xX=\on{Spec}\bigoplus_{k\in\BN}\fm_x^k/\fm_x^{k+1}$. However, the kernel
of the $\fg$-action on $C_xX$ is the Lie algebra of the kernel of the $G$-action
on $C_xX$ (which in turn coincides with the kernel of the $G$-action on $X$,
due to the reductivity of $G$).
The nilpotent ideal $\fq^{>0}$ integrates to
a unipotent algebraic group $\sQ^{>0}$, and the adjoint action of $\fq^0$
on $\fq^{>0}$ integrates to the action of $G$ on $\fq^{>0}$ (factoring through
$G\twoheadrightarrow\sQ^0$). The lemma follows. \qed

%This means that the restriction of
%Poisson bracket on $T^*_{Zar,p}$ is correctly defined and immediately makes it
%a Lie algebra, which will be further denoted by $ \mf Z $. Grading on $A$ is
%multiplicative, so it defines a grading on $T^*_{Zar,p}$ and for compatibility
%with Lie algebra structure we need to shift it by $k$. In other words,
%$\mf Z^{(i)}=A_{i+k}/\mf m_p^2 $. Also let $G_Z$ be algebraic group with Lie
%algebra $\mf Z$.

%We also need the additional action of reductive group $G$ on $M$, satisfying
%several properties. As later we will need to check them for some examples it
%is convenient to put them in the separate list:
%\begin{itemize}
% \item The action is Hamiltonian.
% \item Corresponding map of Lie algebras $\mf g\rightarrow\mf Z$ is an
%injection.
% \item There is a subtorus $T \hookrightarrow G$ of $G$, with compact form
%$T^c$, such that  $p$ is isolated fixed point of $T^c$-action.
%\end{itemize}
%As a consequence of first two properties we get the existence of a nonempty
%reductive Lie subalgebra of $\mf Z$ and this fact will be crucial for the
%proofs later. It is hard to say right now what for do we need the third
%property, but this should become clear in the next subsection.

\subsection{Calculation of Euler obstructions}
\label{vanish}
There is a Zariski open neighbourhood $U$ of $x$ in $X$, and a closed
embedding $U\hookrightarrow Y$ into a {\em smooth} algebraic variety
which induces an isomorphism on the Zariski tangent spaces $T_xX=T_xU\iso
T_xY$. Dually, $T_x^*Y\iso\fq$. The regular part $\fq^{\on{reg}}\subset\fq$
is defined as $\fq\setminus\bigcup_{i\in L,\ X_i\ne\{x\}}\ol{T^*_{X_i\cap U}Y}\cap
T^*_xY$. This is a nonempty open subset of $\fq$. It is clearly independent
of the choices of $U$ and $Y$ above. Moreover, if we choose any closed
embedding $U\hookrightarrow Y'$ giving rise to a surjection
$p:\ T_x^*Y'\twoheadrightarrow\fq$, then the regular part of $T^*_xY'$
(that is
$T^*_xY'\setminus\bigcup_{i\in L,\ X_i\ne\{x\}}\ol{T^*_{X_i\cap U}Y'}\cap T^*_xY'$)
is nothing but $p^{-1}(\fq^{\on{reg}})$.

\begin{lem}
\label{deform}
$\fq^{\on{reg}}\subset\fq$ is invariant with respect to the adjoint action
of $\sQ$ on $\fq$.
\end{lem}

\proof
Let $C_x\ol{X_i}\subset C_xX\subset T_xX=T_xY$ be the normal cones.
Let $C^{\on{sm}}_x\ol{X_i}\subset C_x\ol{X_i}$ be the smooth part.
The deformation to the normal cone construction proves that
$\ol{T^*_{C^{\on{sm}}_x\ol{X_i}}T_xY}\cap T_x^*Y=\ol{T^*_{X_i\cap U}Y}\cap T^*_xY$.
Now the LHS is the cone over the closed subvariety of $\BP\fq=\BP T_x^*Y$
projectively dual to $\BP C_x\ol{X_i}\subset\BP T_xY$. Since
$C_x\ol{X_i}\subset\fq^*$ is a Poisson subvariety, it is $\sQ$-invariant;
hence its projective dual is $\sQ$-invariant as well.
\qed

\begin{thm}
\label{euler=0}
The Euler obstruction $e_{x,\overline{X}_i}=0$ for every $X_i\ne x$.
\end{thm}

\proof
First we consider $e_{x,X}$.
According to~\cite[2.13]{st}, there is a nonzero $q\in\fq^{\on{reg}}$ stabilized
by a maximal torus $T\subset\sQ$. Since all the maximal tori in $\sQ$ are
$\sQ$-conjugate, while $\fq^{\on{reg}}$ is $\sQ$-invariant, we can find another
nonzero element $q'\in\fq^{\on{reg}}$ stabilized by the maximal torus
$\sA\subset G\to\sQ$ (notations of~\ref{setup}). Choose an
$\sA$-equivariant closed embedding $X\hookrightarrow\BA^N$, and lift $q'$
to an $\sA$-invariant linear function $\phi$ on $\BA^N$. Now the argument
of~\cite[proof of~Theorem~2.2]{em} proves $e_{x,X}=0$.

For an arbitrary symplectic leaf $x\ne X_i\subset X$ we just replace $X$
by $\overline{X}_i$ in the above argument: all the assumptions of~\ref{setup}
are clearly satisfied with the same $\sA\subset G,\ \BC^\times$.
\qed

\section{Applications to transversal slices in double affine Grassmannians}
\label{a}

\subsection{Uhlenbeck spaces}
\label{Uhl}
We follow the notations of~\cite{bf}. To begin with, for an almost simple
simply connected complex algebraic group $H$ we denote by $\CU^a_H(\BA^2)$
the Uhlenbeck closure of the moduli space $\on{Bun}_H^a(\BA^2)$ of $H$-bundles
on $\BP^2$ trivialized along $\BP^1$ of second Chern class $a$. It is acted
upon by $H\times\on{Aff}(2)$. Here $H$ acts by the change of trivialization,
while $\on{Aff}(2)=\on{GL}(2)\ltimes\BG_a^2$ (the group of affine motions of
$\BA^2$, i.e. $\on{Aut}(\BP^2,\BP^1)$) acts via the transport of structure.
The normal subgroup $\BG_a^2$ acts on $\CU^a_H(\BA^2)$ freely, and the
quotient is denoted $'\CU^a_H(\BA^2)$: the reduced Uhlenbeck space,
acted upon by $H\times\on{GL}(2)$. Let $\BC^\times\subset\on{GL}(2)$ stand
for the central subgroup, and let $G$ stand for $H\times\on{SL}(2)$.

\begin{prop}
\label{uhl}
$X=\ '\CU^a_H(\BA^2)$ with the action of $G\times\BC^\times$ satisfes
the assumptions of~\ref{setup}.
\end{prop}

\proof
In case $H=\on{SL}(r)$ the Uhlenbeck space $\CU_{\on{SL}(r)}^a(\BA^2)$
was constructed
by Donaldson by Hamiltonian reduction from the representation space of the
2-loop quiver; it is usually denoted $M_0(V,W)$ where $V=\BC^a,\ W=\BC^r$.
This construction provides $\CU_{\on{SL}(r)}^a(\BA^2)$ with a Poisson structure.
For an arbitrary $H$ and a representation $\varrho:\ H\to\on{SL}(r)$ we
have the corresponding embedding
$\varrho_*:\ \CU^a_H(\BA^2)\hookrightarrow\CU_{\on{SL}(r)}^{a'}(\BA^2)$ for
certain $a'$. The Poisson structure on $\CU_{\on{SL}(r)}^{a'}(\BA^2)$
restricted to $\CU^a_H(\BA^2)$ is independent of $\varrho$ upto a (nonzero)
scalar factor. We fix the Poisson structure on $\CU^a_H(\BA^2)$,
say choosing the adjoint representation $\varrho$.

The well known stratification
$\CU_H^a(\BA^2)=
\bigsqcup_{0\leq b\leq a}\on{Bun}^b_H(\BA^2)\times\on{Sym}^{a-b}\BA^2$ admits
a refinement by the {\em diagonal} stratification of
$\on{Sym}^{a-b}\BA^2=\bigsqcup_{\fP\in P(a-b)}\on{Sym}^\fP\BA^2$,
see~\cite[Section~10]{bfg}
(here $P(a-b)$ stands for the set of partitions of $a-b$). The diagonal strata
are nothing but the symplectic leaves of the Poisson structure; in particular,
there are finitely many symplectic leaves. We will denote the stratum
$\on{Bun}^b_H(\BA^2)\times\on{Sym}^\fP\BA^2$ by $\CU_H^{b,\fP}$ for short,
and its image in the reduced Uhlenbeck space $'\CU^a_H(\BA^2)$ will be
denoted by $'\CU_H^{b,\fP}$. Thus
$'\CU^a_H(\BA^2)=\bigsqcup_{0\leq b\leq a}^{\fP\in P(a-b)}\ '\CU_H^{b,\fP}$
is the decomposition into symplectic leaves. Among those, there is a unique
0-dimensional leaf: the one for $b=0,\ \fP=(a)$. This point $x$ is the unique
fixed point for any maximal torus $\sA\subset G$.

In order to compute the weight $w$ of the Poisson structure with respect to
the $\BC^\times$-action we recall the Hamiltonian reduction construction of
$\CU_{\on{SL}(r)}^a(\BA^2)=M_0(V,W)$. Namely, for $M=\on{Hom}(W,V)\oplus
\on{Hom}(V,W)\oplus\on{End}(V,V)\oplus\on{End}(V,V)$ (with a typical element
$(p,q,A,B)$), and $M\supset M':=\{(p,q,A,B):\ AB-BA+pq=0\}$ naturally acted
upon by $\on{GL}(V)$, we have $M_0(V,W)=M'/\!/\on{GL}(V)$.
The action of $\BG_a^2$ on $M_0(V,W)$ comes from the action
$(a,b)\cdot(p,q,A,B)=(p,q,A+a\on{Id}_V,B+b\on{Id}_V)$ on $M$. Thus
$'\CU_{\on{SL}(r)}^a(\BA^2)=M_0(V,W)/\BG_a^2=M''/\!/\on{GL}(V)$ where
$M'\supset M'':=\{(p,q,A,B):\ AB-BA+pq=0,\ \on{Tr}A=\on{Tr}B=0\}$.
The vector space
$M$ has a natural symplectic structure which gives rise to the Poisson
structure on the categorical quotient $M_0(V,W)$. The action of $\BC^\times$
on $M_0(V,W)$ comes from the dilation action $t\cdot(p,q,A,B)=(tp,tq,tA,tB)$
on $M$. Evidently, the symplectic form on $M$ has weight 2 with respect to
this action, so $w=2$. By the construction of the Poisson structure on
$'\CU_H^a(\BA^2)$, it has weight 2 for arbitrary $H$ as well.

It remains to find all the functions on $'\CU_H^a(\BA^2)$ of
$\BC^\times$-weights 1,2. Again we start with $H=\on{SL}(r)$. By the classical
invariant theory, all the $\on{GL}(V)$-invariant functions on $M''$ are
generated by the following ones: (a) matrix elements of $qCp\in\on{End}(W)$
where $C$ is a word in the alphabet $(A,B)$; (b) traces of $C\in\on{End}(V)$
where $C$ is a word in the alphabet $(A,B)$. Evidently, the $\BC^\times$-weight
of (a) is $\on{length}(C)+2$, while $\on{weight}(\on{Tr}C)=\on{length}(C)$.
Note that the only words of length 1 are $A,B$, and their traces vanish by
definition of $M''$. So there are no functions of weight 1.

The functions of weight 2 are spanned by the matrix elements of $qp$, and
$\on{Tr}A^2,\on{Tr}B^2,\on{Tr}AB$. Note that $\on{Tr}qp=\on{Tr}pq=0$.
Clearly, the first group of weight 2 functions is lifted from the moment
map $M''/\!/\on{GL}(V)\to\mathfrak{sl}(W)^*=\mathfrak{sl}(r)^*$,
while the second one
is lifted from the moment map $M''/\!/\on{GL}(V)\to\mathfrak{sl}(2)^*$.
Thus $\mathfrak{sl}(r)\oplus\mathfrak{sl}(2)\twoheadrightarrow
\BC[\ '\CU^a_{\on{SL}(r)}]_2$.

Finally, for arbitrary $H\stackrel{\varrho}{\to}\on{SL}(r)$, the weight $k$
functions on $'\CU^a_H$ are just the restrictions of weight $k$
functions on $'\CU^{a'}_{\on{SL}(r)}$ under the closed embedding
$\varrho_*:\ '\CU^a_H\hookrightarrow\ '\CU^{a'}_{\on{SL}(r)}$.
The commutative diagram of moment maps
$\begin{CD}
'\CU^a_H @>>\varrho_*> \ '\CU^{a'}_{\on{SL}(r)}\\
@VVV @VVV\\
\fg^*=\fh^*\oplus\mathfrak{sl}(2)^* @<\varrho^*\oplus\on{Id}<<
\mathfrak{sl}(r)^*\oplus\mathfrak{sl}(2)^*
\end{CD}$ completes the proof.
\qed

\subsection{Characteristic cycle of $\on{IC}(\ '\CU^a_H)$}
\label{charcyc}
We choose a closed embedding $'\CU^a_H\hookrightarrow Y$ into a smooth
variety $Y$, and view $\on{IC}(\ '\CU^a_H)$ as a perverse sheaf on $Y$.

\begin{cor}
\label{uhle}
The multiplicity of the conormal bundle $T^*_{\ '\CU_H^{b,\fP}}Y$ in the
characteristic cycle $\on{CC}(\on{IC}(\ '\CU^a_H))$ equals the total
dimension of the stalk of $\on{IC}(\ '\CU^a_H)$ on the stratum
$'\CU_H^{b,\fP}$ (computed in~\cite[Theorem~7.10]{bfg}).
\end{cor}

\proof
The argument of~\cite[Section~2.5]{em} reduces the proof to the vanishing
of the Euler obstructions $e_{'\CU^{b',\fP'}_H,\ '\ol\CU{}^{b,\fP}_H}$ where
$'\ol\CU{}^{b,\fP}_H\supset\ '\CU^{b',\fP'}_H$
stands for the closure of the stratum $'\CU^{b,\fP}_H$.
We first treat the smallest stratum $x=\ '\CU_H^{b',\fP'}$ for $b'=0,\ \fP'=(a)$.
This vanishing follows
from~Theorem~\ref{euler=0} and~Proposition~\ref{uhl}.
In general, the vanishing $e_{'\CU^{b',\fP'}_H,\ '\ol\CU{}^{b,\fP}_H}=0$ is equivalent
to the vanishing $e_{\CU^{b',\fP'}_H,\ol\CU{}^{b,\fP}_H}=0$ where
$\ol\CU{}^{b,\fP}_H\supset\CU^{b',\fP'}_H$
stands for the closure of the stratum $\CU^{b,\fP}_H$ in the nonreduced
Uhlenbeck space $\CU^a_H(\BA^2)$. By the factorization
principle~\cite[Proposition~6.5]{bfg}, the desired obstruction is the
product $\prod_{l=1}^me_{\CU^{0,(a_l)}_H,\ol\CU{}^{b_l,\fP_l}_H}$ (where $m$ is the number
of parts of the partition $\fP$). The latter factors
are already proved to vanish.   \qed

%Equivalently, we have
%to prove the vanishing of the Euler obstructions in the nonreduced Uhlenbeck
%space: $e_{\CU^{0,(a)}_H,\ol\CU{}^{b,\fP}_H}=0$  where $\ol\CU{}^{b,\fP}_H$
%stands for the closure of the stratum $\CU^{b,\fP}_H$. If
%$\ol{\on{Sym}}{}^\fP\BA^2$
%denotes the closure of a diagonal stratum in $\on{Sym}\BA^2$, then we have
%a finite, generically one-to-one, morphism
%$\CU^b_H(\BA^2)\times\ol{\on{Sym}}{}^\fP\BA^2\to\ol\CU{}^{b,\fP}_H$
%(it is also bijective over the smallest stratum). The preimage of the
%smallest stratum $\CU^{0,(a)}_H\subset\ol\CU{}^{b,\fP}_H$ in
%$\CU^b_H(\BA^2)\times\ol{\on{Sym}}{}^\fP\BA^2$ is the diagonal
%$\Delta\subset\BA^2\times\BA^2=\CU^{0,(b)}_H\times\on{Sym}^{(b-a)}\BA^2\subset
%\CU^b_H(\BA^2)\times\ol{\on{Sym}}{}^\fP\BA^2$. It suffices to prove the
%vanishing $e_{\Delta,\CU^b_H(\BA^2)\times\ol{\on{Sym}}{}^\fP\BA^2}=0$.
%Taking quotient by the free action of $\BG_a^2$ once again

\subsection{Transversal slices in double affine Grassmannians}
Given a cyclic subgroup $\Gamma_k=\BZ/k\BZ\subset G=H\times\on{SL}(2)$,
a transversal slice $\CU^\lambda_{H,\mu}(\BA^2/\Gamma_k)$ is defined
in~\cite{bf} as a certain irreducible component of the fixed point set
$\CU^a_H(\BA^2)^{\Gamma_k}$. The Poisson structure is restricted from the one on
$\CU^a_H(\BA^2)$, and the symplectic leaves are the intersections of
$\CU^\lambda_{H,\mu}(\BA^2/\Gamma_k)$ with the symplectic leaves (the strata of
the diagonal stratification) of $\CU^a_H(\BA^2)$.

In case $k=2$, the cyclic subgroup $\Gamma_2$ is central in $\on{SL}(2)$,
and the centralizer
$G=Z_{H\times\on{SL}(2)}(\Gamma_2)=Z_H(\Gamma_2)\times\on{SL}(2)$ acts on
$\CU^\lambda_{H,\mu}(\BA^2/\Gamma_2)$ preserving the Poisson structure.
In case $k>2$, we only have the action of $G=Z_{H\times\on{SL}(2)}(\Gamma_k)=
Z_H(\Gamma_k)\times T$ on
$\CU^\lambda_{H,\mu}(\BA^2/\Gamma_k)$ where $T$ is the centralizer torus of
$\Gamma_k$ in $\on{SL}(2)$.
The action of central $\BC^\times\subset\on{GL}(2)$
on $\CU^a_H(\BA^2)$ preserves $\CU^\lambda_{H,\mu}(\BA^2/\Gamma_k)$.

\begin{prop}
\label{slice}
$X=\CU^\lambda_{H,\mu}(\BA^2/\Gamma_k)$ with the action of $G\times\BC^\times$
satisfies the assumptions of~\ref{setup}.
\end{prop}

\proof
We essentially repeat the proof of~Proposition~\ref{uhl}: using a
representation $\varrho:\ H\to\on{SL}(r)$ we reduce the claim to the case
$H=\on{SL}(r)$. In this case the transversal slice
$\CU^\lambda_{H,\mu}(\BA^2/\Gamma_k)$ is nothing but a cyclic quiver variety
$M_0(\ul{V},\ul{W})$~\cite[Section~7]{bf} (for the cyclic quiver with $k$
vertices). We have $M_0(\ul{V},\ul{W})=M'/\!/\on{GL}(\ul{V})$ where
$\on{GL}(\ul{V})=\prod_{l\in\BZ/k\BZ}\on{GL}(V_l)$, and $M'\subset M:=
\bigoplus_{l\in\BZ/k\BZ}\on{Hom}(W_l,V_l)\oplus\bigoplus_{l\in\BZ/k\BZ}
\on{Hom}(V_l,W_l)\oplus\bigoplus_{l\in\BZ/k\BZ}\on{End}(V_l,V_{l+1})\oplus
\bigoplus_{l\in\BZ/k\BZ}\on{End}(V_l,V_{l-1})$ is cut out by the equations
$A_{l-1}B_l-B_{l+1}A_l+p_lq_l=0,\ l\in\BZ/k\BZ$. The $\BC^\times$-action is by
dilations, and the $\BC^\times$-weight of the Poisson structure is $w=2$.
It is not hard to check from the classical invariant theory that
all the $\on{GL}(\ul{V})$-invariant functions on $M'$ are generated by the
following ones: (a) matrix elements of $q_nCp_m$ where $C$ is a word in
the alphabet $(A_l,B_l)_{l\in\BZ/k\BZ}$ (not all the words are allowed: only the
{\em composable} ones); (b) traces of $C\in\on{End}(V_0)$ where $C$ is
a composable word in the alphabet $(A_l,B_l)_{l\in\BZ/k\BZ}$ starting and ending
at the 0-th vertex. Among those, the functions of weight 2 are the matrix
elements of $q_lp_l$, and $\on{Tr}(B_1A_0)$ for $k>2$ (note that
$\on{Tr}(A_{-1}B_0)=\on{Tr}(B_1A_0)-\on{Tr}(q_0p_0$), and also $\on{Tr}A_1A_0,\
\on{Tr}B_1B_0$ for $k=2$. Clearly, these functions are lifted from the
moment map $M_0(\ul{V},\ul{W})\to\left(\bigoplus_{l\in\BZ/k\BZ}\mathfrak{gl}(W_l)
\cap\mathfrak{sl}(\bigoplus_{l\in\BZ/k\BZ}W_l)\right)^*
\oplus\ft^*$ in case $k>2$, and
$M_0(\ul{V},\ul{W})\to\left(\bigoplus_{l\in\BZ/k\BZ}\mathfrak{gl}(W_l)
\cap\mathfrak{sl}(\bigoplus_{l\in\BZ/k\BZ}W_l)\right)^*
\oplus\mathfrak{sl}(2)^*$ in case $k=2$.
The proposition is proved. \qed

\begin{cor}
\label{slic}
For a closed embedding $\CU^\lambda_{H,\mu}(\BA^2/\Gamma_k)\hookrightarrow Y$
into a smooth variety $Y$, the multiplicity of the conormal bundle
$T^*_SY$ to a symplectic leaf $S\subset\CU^\lambda_{H,\mu}(\BA^2/\Gamma_k)$
in the characteristic cycle $\on{CC}(\on{IC}(\CU^\lambda_{H,\mu}(\BA^2/\Gamma_k)))$
equals the total dimension of the stalk of
$\on{IC}(\CU^\lambda_{H,\mu}(\BA^2/\Gamma_k))$ on the stratum $S$.
\end{cor}

\proof The same as the proof of~Corollary~\ref{uhle}. \qed

\begin{rem}
{\em The dimension of the stalk of $\on{IC}(\CU^\lambda_{H,\mu}(\BA^2/\Gamma_k))$
on a stratum $S$ is computed by the factorization principle,
and~\cite[Conjecture~4.14]{bf}.}
\end{rem}

\section{A conjecture}
\label{aco}

\subsection{Chambers}
\label{chambers}
In the setup of~\ref{setup} we set $\fa:=\on{Lie}\sA$, and
$\fa_\BR:=X_*(\sA)\otimes_\BZ\BR\subset\fa$. We say that a coweight
$a\in X_*(\sA)$ is {\em regular} if $X^{a(\BC^\times)}=x$. We assume that the
nonregular coweights form a finite union of corank one subgroups in
$X_*(\sA)$. We say that $a\in\fa_\BR$ is regular if it lies off the
corresponding real hyperplanes. The connected components of $\fa_\BR^{\on{reg}}$
are called {\em chambers}.

Given a regular $a\in X_*(\sA)$ we define the attracting set
$\fT_a\subset X$ as the set of all $z\in X$ such that
$\lim_{c\to0}a(c)\cdot z=x$.
If $a$ varies in a chamber $\fC$, then $\fT_a$ does not change, and so we
denote it $\fT_\fC$. The closed embedding $\fT_\fC\hookrightarrow X$ is
denoted by $\iota_\fC$. The closed embedding $x\hookrightarrow\fT_\fC$ is denoted
by $j_\fC$. The composition $j_\fC^!\iota_\fC^*$ is the hyperbolic restriction.
In all the examples of~Section~\ref{a} $j_\fC^!\iota_\fC^*\on{IC}(X)$ is a
vector space in cohomological degree 0.

\subsection{Microlocalization}
\label{micro}
According to the first paragraph of~Section~\ref{vanish}, we have a
perverse sheaf $\bmu(\on{IC}(X))$ on $T^*_x=\fq$. According
to~\cite[Proposition~4.4.7]{ks}, $\bmu(\on{IC}(X))$ is well defined, i.e.
does not depend on the choice of $U\hookrightarrow Y$.
Recall that we have a homomorphism $\eta:\ \fa\to\fq=T_x^*X$.
Motivated by~\cite{mo} and by the examples of~Section~\ref{a} we propose the
following

\begin{conj}
\label{oko}
(a) $\eta^*\bmu(\on{IC}(X))$ is constant in any chamber $\fC$.

(b) There is a canonical isomorphism $\eta^*\bmu(\on{IC}(X))_\fC\iso
j_\fC^!\iota_\fC^*\on{IC}(X)$.
\end{conj}

\end{document}